\pdfoutput=1
\RequirePackage{ifpdf}
\ifpdf 
\documentclass[pdftex]{sigma}
\else
\documentclass{sigma}
\fi

\numberwithin{equation}{section}
\newtheorem{Theorem}{Theorem}[section]
\newtheorem{Corollary}[Theorem]{Corollary}
\newtheorem{Lemma}[Theorem]{Lemma}
\newtheorem{Proposition}[Theorem]{Proposition}
{ \theoremstyle{definition}
\newtheorem{Definition}[Theorem]{Definition}
\newtheorem{Notation}[Theorem]{Notation}
\newtheorem{Example}[Theorem]{Example}
\newtheorem{Remark}[Theorem]{Remark} }

\usepackage{stmaryrd}
\usepackage[all]{xy}
\usepackage{tikz}\usetikzlibrary{arrows}
\DeclareMathAlphabet{\mathpzc}{OT1}{pzc}{m}{it}
\usepackage{mathrsfs}
\newcommand{\op}[1]{{#1}^{\mbox{\sf{\tiny{op}}}}}
\newcommand{\gf}[1]{{#1}_{\mbox{\sf{\tiny{f\/ib}}}}}
\newcommand{\opfib}[1]{{{#1}_{\mbox{\sf{\tiny{f\/ib}}}}}^{\mbox{\sf{\tiny{op}}}}}

\begin{document}

\newcommand{\arXivNumber}{1208.1513}

\allowdisplaybreaks

\renewcommand{\PaperNumber}{022}

\FirstPageHeading

\ShortArticleName{Dynamics on Networks of Manifolds}

\ArticleName{Dynamics on Networks of Manifolds}

\Author{Lee DEVILLE and Eugene LERMAN}

\AuthorNameForHeading{L.~DeVille and E.~Lerman}
\Address{Department of Mathematics, University of Illinois, USA}
\Email{\href{mailto:rdeville@illinois.edu}{rdeville@illinois.edu},
\href{mailto:lerman@illinois.edu}{lerman@illinois.edu}}
\URLaddress{\url{http://www.math.illinois.edu/~rdeville/},\\
\hspace*{10.5mm}\url{http://www.math.illinois.edu/~lerman/}}

\ArticleDates{Received April 24, 2014, in f\/inal form February 24, 2015; Published online March 12, 2015}

\Abstract{We propose a~precise def\/inition of a~continuous time dynamical system made up of interacting open subsystems.
The interconnections of subsystems are coded by directed graphs.
We prove that the appropriate maps of graphs called {\it graph fibrations} give rise to maps of dynamical systems.
Consequently surjective graph f\/ibrations give rise to invariant subsystems and injective graph f\/ibrations give rise to
projections of dynamical systems.}


\Classification{34C14; 18D99}

\section{Introduction}

Given a~dynamical system, one often starts by trying to f\/ind invariant subsystems; these include equilibria, periodic
orbits, and higher dimensional invariant submanifolds.
In addition, constructing projections onto smaller systems as well as conjugacies and semi-conjugacies with simpler
systems are generally useful for understanding the qualitative properties of dynamical systems.
All of these objects: invariant subsystems, projections, conjugacies and semi-conjugacies can be realized as maps of
dynamical systems (q.v.~Def\/inition~\ref{def:dyn-sys}).
Thus the search for maps between dynamical systems may be considered one of the fundamental questions of the subject.

In this paper we give a~precise def\/inition of a~continuous time dynamical system made up of interacting open subsystems.
We then exploit the combinatorial aspect of such systems to produce maps of dynamical systems out of appropriate maps of
graphs called {\it graph fibrations} (q.v.\ Def\/inition~\ref{graph fibration}).
We show that in particular surjective graph f\/ibrations give rise to invariant subsystems and injective graph f\/ibrations
give rise to projections of dynamical systems.

The present work is part of an ongoing project.
In~\cite{DL} we reformulated the groupoid formalism of Golubitsky, Pivato, Stewart and
T{\"o}r{\"o}k~\cite{Golubitsky.Stewart.Torok.05, Stewart.Golubitsky.Pivato.03} for coupled cell networks (which are
systems of ordinary dif\/ferential equations) in a~coordinate free manner and extended it to groupoid-invariant vector f\/ields on manifolds.
A~preliminary version was posted as~\cite{DLold}.
We later realized that groupoid invariance of vector f\/ields is not needed for the existence of invariant subspaces.
With the benef\/it of hindsight we see that the theory developed in~\cite{DL} is an equivariant version of the theory that
we develop here.
We would like to point out that dropping groupoid invariance makes the theory much simpler and more f\/lexible.
In particular we expect the results of this paper to readily generalize to hybrid systems.

The absence of explicit groupoid symmetries makes our work close in spirit to the approach to dynamics on networks
advocated by Field~\cite{FieldCD}.
Unlike Field we f\/ind it convenient to use the language of category theory.
We also f\/ind it useful to borrow the notions of open systems and their interconnection from engineering (see, for
example~\cite{Brockett, Pappas, willems}) and the def\/inition of a~graph f\/ibration from computer science~\cite{Vigna1}
(see~\cite{VignaFibPage} for a~history of the notion and alternative terminologies).

We believe that both Field's approach and ours is based on the existence of a~certain algebraic structure which at the
present time is not completely understood.
Open continuous time systems form an algebra over a~certain operad~\cite{VSL}.
This operad is implicit in the work of Field~\cite{FieldCD}.
A~piece of this algebra shows up in our work as the interconnection maps (see Theorem~\ref{thm:3.8}).
We do not understand yet how graph f\/ibrations interact with this operad and plan to address this issue in a~future work.

The goal of this paper is to construct a~category of networks of continuous time systems and a~functor to the category
of dynamical systems.
A~network in our sense consists of
\begin{itemize}\itemsep=0pt

\item a~f\/inite directed graph~$G$ with a~set of nodes $G_0$,

\item a~{\em phase space function} ${\mathcal P}$ that assigns to each node of the graph an appropriate phase space
(which we take to be a~manifold),

\item a~family of open systems $\{w_a\}_{a\in G_0}$ (one for each node~$a$ of the graph~$G$) consistent in an
appropriate way with the structure of the graph, and

\item an interconnection map $\mathscr{I}$ that turns these open systems into a~vector f\/ield on the product
$\bigsqcap\limits_{a\in G_0} {\mathcal P} (a)$ of the phase spaces of the nodes.
\end{itemize}
Our main result, Theorem~\ref{thm:main}, shows that graph f\/ibrations compatible with phase space functions give rise to
maps of dynamical systems.
This allows us to def\/ine a~functor from dynamical systems on networks to general dynamical systems.

The reader may wonder what motivates us to come up with these def\/initions and constructions.
Indeed there are many dif\/ferent kinds of objects in engineering, science and mathematics that are called ``networks''.
The notion of a~network in the present paper arose from the follo\-wing idea, which is implicit in the literature on
coupled cell networks.
Imagine a~physical system modeled by a~vector f\/ield~$X$ on a~manifold~$M$;~$M$ is the collection of all possible states
of the system.
Such systems are common in classical mechanics, to give one example.
Suppose further that our system consists of two interacting subsystems.
We can model this by saying that the collection of states of the f\/irst subsystem forms a~manifold $M_1$ and the second
a~manifold $M_2$.
We would like the states of the big system to be completely determined by the states of its subsystems.
We model this by requiring that $M = M_1 \times M_2$.
A~vector f\/ield~$X$ on $M_1\times M_2$ then has to be of the form
\begin{gather*}
X(x_1,x_2) = (X_1(x_1,x_2), X_2(x_1,x_2)),
\end{gather*}
where
\begin{gather*}
X_1(x_1,x_2) \in T_{x_1}M_1
\qquad
\text{for all}
\quad
(x_1,x_2)\in M_1\times M_2,
\end{gather*}
with a~similar equation holding for $X_2\colon M_1 \times M_2 \to TM_2$.
Note that the functions $X_1$, $X_2$ are {\em not} vector f\/ields.
They are open systems in the sense of Def\/inition~\ref{def:open-sys}.
Moreover the vector f\/ield~$X$ may be considered to be a~result of interconnecting $X_1$ and $X_2$ (see
Proposition~\ref{prop:3.2} and Theorem~\ref{thm:3.8}).

To continue with our example, observe that the evolution of the subsystem 1 depends on its state and the state of the subsystem~2.
Similarly the evolution the second subsystem depends on its state and the state of the subsystem 1.
These mutual inf\/luences can be pictured graphically~as
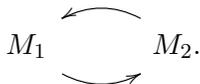
\begin{gather*}
\xy
(-10,0)*++{M_1}="1";
(10,0)*++{M_2.}="2";
{\ar@/_1.2pc/ "1";"2"};
{\ar@/_1.2pc/ "2";"1"};
\endxy
\end{gather*}

Assume now that the map $X_2$ {\em does not} really depend on the points of $M_1$.
That is, there is a~map $Y:M_2 \to TM_2$ with $Y(x_2) \in T_{x_2}M_2$ and $X_2(x_1,x_2) = Y(x_2)$ for all $(x_1,x_2)\in
M_1\times M_2$.
We can picture this as
\begin{gather*}
\xy
(-10,0)*++{M_1}="1";
(10,0)*++{M_2}="2";
{\ar@{->} "2";"1"};
\endxy
\end{gather*}
and say that the second subsystem drives the f\/irst but not conversely.
This way of picturing a~system made up of interacting subsystems generalizes to any number of subsystems.
For example, a~system may be made up of three interacting subsystems like this:
\begin{gather}
\label{eq:favorite}
\xy
(-15,0)*++{M_1}="1";
(15,0)*++{M_2}="2";
(45,0)*++{M_3.}="3";
{\ar@/_1.5pc/ "1";"2"};
{\ar@{->}_{} "2";"3"};
{\ar@/_1.5pc/ "2";"1"};
\endxy
\end{gather}
The total phase space of such a~system would be the product $M=M_1\times M_2 \times M_3$ and the dynamics would be
governed by a~vector f\/ield~$X$ of the form
\begin{gather*}
X(x_1,x_2, x_3) = (X_1(x_1, x_2), X_2 (x_2, x_1), X_3(x_3, x_2)).
\end{gather*}
How are we then to interpret the diagram of the form $ M \xymatrix{*+{ } \ar@(dr,ur)}$~~~~?
And why would we want to? Here is a~two part answer.
We interpret this diagram as a~vector f\/ield~$X$ on the manifold~$M$ of the form
\begin{gather*}
X(x) = w(x,x),
\end{gather*}
where $w:M\times M\to TM$ is an open system with $w(x_1,x_2) \in T_{x_1} M$ for all $(x_1,x_2) \in M\times M$.
This seems a~bit strange and pedantic, but it is useful.
Consider a~vector f\/ield~$Z$ on $M\times M\times M$ of the form
\begin{gather*}
Z(x_1,x_2, x_3) = (w(x_1, x_2), w(x_2,x_1), w(x_3,x_2)),
\end{gather*}
where~$w$ is the open system above.
The vector f\/ield~$Z$ on~$M$ models the dynamics of a~system consisting of three interacting subsystems with the f\/irst
driving the second, the second driving the f\/irst and the third just as in~\eqref{eq:favorite}, only now all the
subsystems have isomorphic phase spaces.
It is not hard to check that the diagonal
\begin{gather*}
\Delta_M:= \{(x_1,x_2,x_3)\in M\times M\times M \,|\, x_1 = x_2 = x_3\}
\end{gather*}
is an invariant submanifold for the vector f\/ield~$Z$.
According to the philosophy we brought up in the f\/irst paragraph of the paper the invariance of $\Delta_M$ should be
seen as coming from a~map of dynamical systems.
And indeed the diagonal map
\begin{gather*}
\delta \colon \ M\to M\times M\times M,
\qquad
\delta(x) = (x,x,x)
\end{gather*}
gives rise to a~map of dynamical systems $\delta\colon (M, X) \to (M\times M\times M, Z)$.
The main result of the paper, Theorem~\ref{thm:main}, implies that this map of dynamical systems is induced by the map
of graphs
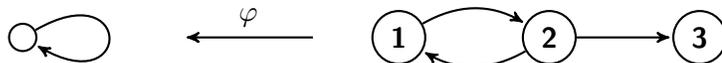
\begin{gather*}
\tikzset{every loop/.style={min distance=15mm,in=-30,out=30,looseness=10}}
\begin{tikzpicture}
[->,>=stealth',shorten >=1pt,auto,node distance=2cm, thick,main node/.style={circle,draw,font=\sffamily\bfseries}]
\node[main node] at (2,-1) (1) {1}; \node[main node] at (4,-1) (2) {2}; \node[main node] at (6,-1) (3) {3}; \node at
(1,-1) (rb){}; \node at (-1,-1) (lb){}; \node[main node] at (-3,-1) (circ){};
\path[] (1) edge [bend left] node {} (2) (2) edge [bend left] node [below] {} (1) (2) edge node {}(3) (circ) edge [loop
right=20] node {} (circ) (rb) edge node [above] {$\varphi$} (lb);
\end{tikzpicture}
\end{gather*}
which is a~graph f\/ibration.
We note that the vector f\/ield~$Z$ has groupoid symmetry in the sense of Golubitsky et al.~\cite{Golubitsky.Stewart.Torok.05,Stewart.Golubitsky.Pivato.03} and~\cite{DL}.
For us, however, the groupoid invariance of~$Z$ is, in some sense, incidental.
It is a~consequence of the fact that~$Z$ is assembled out of the triple of open systems which lies in the image of the
map $\varphi^*$ of Theorem~\ref{thm:4.4} and that~$\varphi$ happens to be surjective.

The paper is organized as follows.
We start by def\/ining the category $\mathsf{DS}$ of continuous time dynamical systems.
We recall the def\/inition of a~directed multigraph, def\/ine the notion of a~network of manifolds and the total space of
the network.
We recall the notion of an open system and discuss interconnections of open systems.
We show how a~network of manifolds naturally leads to a~collection of spaces of open systems that can be interconnected.
We then prove our main result, Theorem~\ref{thm:main}: f\/ibrations of networks of manifolds give rise to maps of
dynamical systems.
We end the paper with a~collection of examples.

\section{Definitions and constructions}

\noindent We start by def\/ining what we mean by a~continuous time dynamical system and by a~map between two such systems.
\begin{Definition}
\label{def:dyn-sys}
A~continuous time {\it dynamical system} is a~vector f\/ield on a~manifold.
More formally it is a~pair $(M,X)$, where~$X$ is a~vector f\/ield on a~manifold~$M$.

A {\it map} from a~dynamical system $(M,X)$ to a~dynamical system $(N, Y)$ is a~smooth map $f:M\to N$ that intertwines
the two vector f\/ields:
\begin{gather*}
Df\circ X = Y\circ f,
\end{gather*}
where $Df:TM\to TN$ denotes the dif\/ferential of~$f$.
One also says that the vector f\/ields~$X$ and~$Y$ are~$f$-related.
\end{Definition}

\begin{Notation}[the category $\mathsf{DS}$ of dynamical systems]
\label{rmrk:dyn-sys}
Continuous time dynamical systems and maps of dynamical systems form a~category.
We denote it by $\mathsf{DS}$.
\end{Notation}

\subsection{Graphs and manifolds}

Throughout the paper {\em graphs} are f\/inite directed multigraphs, possibly with loops.
More precisely, we use the following def\/inition:
\begin{Definition}
A~{\it graph} $ G$ consists of two f\/inite sets $G_1$ (of arrows, or edges), $G_0$ (of nodes, or vertices) and two maps
$\mathfrak s,\mathfrak t\colon G_1 \to G_0$ (source, target); we write
\begin{gather*}
G = \{G_1\rightrightarrows G_0\}.
\end{gather*}
The set $G_1$ may be empty, i.e., we may have $G= \{\varnothing \rightrightarrows G_0\}$, making~$G$ a~disjoint
collection of vertices with no arrows between them.
\end{Definition}

\begin{Definition}
A~{\it map of graphs} $\varphi\colon A\to B$ from a~graph~$A$ to a~graph~$B$ is a~pair of maps $\varphi_1\colon A_1\to
B_1$, $\varphi_0\colon A_0\to B_0$ taking edges of~$A$ to edges of~$B$, nodes of~$A$ to nodes of~$B$ so that for any
edge~$\gamma$ of~$A$ we have
\begin{gather*}
\varphi_0 (\mathfrak s(\gamma)) = \mathfrak s(\varphi_1 (\gamma))
\qquad
\text{and}
\qquad
\varphi_0 (\mathfrak t(\gamma)) = \mathfrak t(\varphi_1 (\gamma)).
\end{gather*}
We often omit the indices 0 and 1 and write $\varphi(\gamma)$ for $\varphi_1 (\gamma)$ and $\varphi(a)$ for $\varphi_0(a)$.
\end{Definition}

\begin{Remark}
The collection of f\/inite (directed multi-)graphs and maps of graphs form a~category $\mathsf{Graph}$.
\end{Remark}

In order to construct networks from graphs we need to have a~consistent way of assigning manifolds to nodes of our
graphs.
We formalize this idea by making the collection of graphs with manifolds assigned to vertices into a~category
$\mathsf{Graph}/\mathsf{Man}$.

\begin{Definition}[category of networks of manifolds $\mathsf{Graph}/\mathsf{Man}$]\label{def.network.mfld}
A~{\it network of manifolds} is a~pair $(G,{\mathcal P})$, where~$G$
is a~(f\/inite directed multi-)graph and ${\mathcal P}\colon G_0 \to \mathsf{Man}$ is a~function that assigns to each node~$a$ of~$G$
a~manifold ${\mathcal P}(a)$.
We think of ${\mathcal P}$ as an assignment of phase spaces to the nodes of the graph~$G$, and for this reason we refer
to ${\mathcal P}$ as a~{\em phase space function.}

Networks of manifolds form a~category $\mathsf{Graph}/\mathsf{Man}$.
Its objects are are pairs $(G,{\mathcal P})$ as above.
A~morphism~$\varphi$ from $(G,{\mathcal P})$ to$(G',{\mathcal P}')$ is a~map of graphs $\varphi\colon G\to G'$ with
\begin{gather*}
{\mathcal P}'\circ \varphi = {\mathcal P}.
\end{gather*}
\end{Definition}

\begin{Notation}
Given a~category $\mathscr{C}$ we denote the opposite category by $\op{\mathscr{C}}$, i.e.~the category with all of the
same objects and all of the arrows reversed.
We adhere to the convention that a~{\em contravariant functor} from a~category $\mathscr{C}$ to a~category $\mathscr{D}$
is a~covariant functor
\begin{gather*}
F\colon \ \op{\mathscr{C}}\to \mathscr{D}.
\end{gather*}
Then for any morphism $c\xrightarrow{\gamma}c'$ of $\mathscr{C}$ we have $F(c)\xleftarrow{F(\gamma)}F(c')$ in
$\mathscr{D}$.
\end{Notation}

Next we recall the notion of a~{\em product} in a~category $\mathscr{C}$.
We will use them in two instances: when $\mathscr{C}$ is the category $\mathsf{Man}$ of smooth f\/inite dimensional
manifolds and smooth maps and when~$\mathscr{C}$ is the category $\mathsf{Vect}$ of real (but not necessarily f\/inite
dimensional) vector spaces and linear maps.

\begin{Definition}
A~{\it product} of a~family $\{c_s\}_{s\in S}$ of objects in a~category $\mathscr{C}$ indexed by a~set~$S$ is an object
$\bigsqcap_{s'\in S}c_{s'}$ of $\mathscr{C}$ together with a~family of morphisms $\{\pi_s\colon \bigsqcap_{s'\in
S}c_{s'} \to c_s\}_{s\in S}$ with the following universal property: given an object $c'$ of $\mathscr{C}$ and a~family
of morphisms $\{f_s\colon c'\to c_s\}_{s\in S}$ there is a~unique morphism $f\colon c'\to \bigsqcap_{s\in S}c_s$ with
\begin{gather*}
\pi_s \circ f = f_s
\qquad
\text{for all}
\quad
s\in S.
\end{gather*}
\end{Definition}

\begin{Remark}
If a~product exists then it is unique up to a~unique isomorphism~\cite{Awodey}.
\end{Remark}

\begin{Lemma}
The category of manifolds $\mathsf{Man}$ has $($finite$)$ categorical products.
\label{empt:cat-prod}
\end{Lemma}

\begin{proof}
There are several ways to {\em construct} categorical products in $\mathsf{Man}$.
The f\/irst one uses Cartesian products: given a~family $\{M_s\}_{s\in S}$ of manifolds indexed by an~$n$-element set~$S$,
order the elements of~$S$: $S=\{s_1, \ldots, s_n\}$.
Set
\begin{gather*}
\bigsqcap_{s\in S}M_{s} = \prod\limits_{i=1}^n M_{s_i},
\end{gather*}
where the right hand side is the Cartesian product.
The projections $p_{s_j}: \prod\limits_{i=1}^n M_{s_i} \to M_{s_j} $ are just projections on the~$j$-th factor.
It is easy to check that a~product constructed this way has the requisite universal property.
In particular, if we choose two dif\/ferent orderings of elements of~$S$, the resulting products are canonically
isomorphic.
This construction is convenient for writing down examples.

However, for proving the results below, such as Proposition~\ref{prop2.13}, it is better to have a~construction of the
product that does not involve a~choice of ordering of the indexing set in question.
This may be done as follows.
Given a~family $\{M_s\}_{s\in S} $ of manifolds, denote by $\bigsqcup_{s\in S}M_s$ their disjoint union\footnote{It may
be def\/ined by $\bigsqcup\limits_{s\in S}M_s = \bigcup\limits_{s\in S} M_s\times \{s\}$.}.
Now def\/ine
\begin{gather*}
\bigsqcap_{s\in S}M_{s}:= \bigg\{x\colon S \to \bigsqcup_{s\in S}M_s \, \bigg| \,  x(s)\in M_s\ \text{for all}\ s\in S\bigg\}.
\end{gather*}
The projection maps $\pi_s\colon \bigsqcap\limits_{s'\in S}M_{s'}\to M_s$ are def\/ined~by
\begin{gather*}
\pi_s (x) = x(s).
\end{gather*}
We denote $x(s)\in M_s$ by $x_s$ and think of it as the $s^{\rm th}$ ``coordinate''
of an element $x\in \bigsqcap\limits_{s\in S}M_{s}$.
Equivalently we may think of elements of the categorical product $\bigsqcap_{s\in S}M_{s} $ as {\it unordered} tuples
$(x_s)_{s\in S}$ with $x_s\in M_s$.
\end{proof}

\begin{Lemma}
The category of vector spaces $\mathsf{Vect}$ has finite categorical products.
\label{empt:cat-prod-vect}
\end{Lemma}
\begin{proof}[Sketch of proof] Just as in the proof of Lemma~\ref{empt:cat-prod} the f\/inite products in $\mathsf{Vect}$ can be
constructed as vector spaces of ordered tuples of vectors, that is, as Cartesian products.
Categorical products in $\mathsf{Vect}$ can also be constructed as {\em unordered} tuples of vectors.
\end{proof}

\begin{Definition}[total phase space of a~network $(G,{\mathcal P})$] For a~pair $(G,{\mathcal P})$ consisting of a~f\/inite graph~$G$ and
an assignment ${\mathcal P}\colon G_0\to \mathsf{Man}$, that is, for an object $(G,{\mathcal P})$ of
$\mathsf{Graph}/\mathsf{Man}$ we set
\begin{gather*}
\mathbb{P} G\equiv \mathbb{P} (G,{\mathcal P}):= \bigsqcap_{a\in G_0} {\mathcal P} (a),
\end{gather*}
the categorical product of manifolds attached to the nodes of the graph~$G$ by the {\it phase space function} ${\mathcal
P}$ and call the resulting manifold $\mathbb{P} G$ the {\it total phase space} of the network $(G,{\mathcal P})$.
\end{Definition}

\begin{Example}
\label{example:1}
Consider the graph

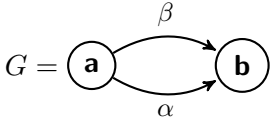
\begin{tikzpicture}
[->,>=stealth',shorten >=1pt,auto,node distance=2cm, thick,main node/.style={circle,draw,font=\sffamily\bfseries}]

\node at (0,0) {$G = $}; \node[main node] at (0.8,0) (1) {a}; \node[main node] at (2.8,0) (2) {b};

\path[every node/.style={font=\sffamily\small}] (1) edge [bend left] node {$\beta$} (2) (1) edge [bend right] node
[below] {$\alpha$} (2);
\end{tikzpicture}

Def\/ine ${\mathcal P}\colon G_0\to \mathsf{Man}$ by ${\mathcal P}(a) = S^2$ (the two sphere) and ${\mathcal P}(b) = S^3$.
Then
\begin{gather*}
\mathbb{P} (G, {\mathcal P}) = S^2 \times S^3.
\end{gather*}
\end{Example}

\begin{Notation}
\label{notation:2.11}
If $G=\{\varnothing \rightrightarrows \{a\}\}$ is a~graph with one node~$a$ and no arrows, we write $G=\{a\}$.
Then for any phase space function ${\mathcal P}\colon G_0 =\{a\} \to \mathsf{Man}$ we abbreviate $\mathbb{P}
(\{\varnothing \rightrightarrows \{a\} \}, {\mathcal P}\colon \{a\} \to \mathsf{Man})$ as $\mathbb{P} a$.
\end{Notation}

\begin{Proposition}
\label{prop2.13}
The assignment
\begin{gather*}
(G,{\mathcal P}) \mapsto \mathbb{P} G:= \bigsqcap_{a\in G_0} {\mathcal P} (a)
\end{gather*}
of phase spaces to networks extends to a~contravariant functor
\begin{gather*}
\mathbb{P} \colon \ \op{(\mathsf{Graph}/\mathsf{Man})} \to \mathsf{Man}.
\end{gather*}
\end{Proposition}

\begin{proof}
Suppose $\varphi\colon (G,{\mathcal P})\to (G',{\mathcal P}')$ is a~morphism in $\mathsf{Graph}/\mathsf{Man}$.
That is, suppose $\varphi\colon G\to G'$ is a~map of graphs with ${\mathcal P}'\circ \varphi = {\mathcal P}$.
We need to def\/ine a~map of manifolds
\begin{gather*}
\mathbb{P} \varphi\colon \ \mathbb{P} G' \to \mathbb{P} G.
\end{gather*}
Since by def\/inition $\mathbb{P} G$ is the product $\bigsqcap\limits_{a\in G_0} {\mathcal P}(a)$, the universal property of
products implies that in order to def\/ine $\mathbb{P} \varphi$ it is enough to def\/ine a~family of maps
\begin{gather*}
\left\{(\mathbb{P} \varphi)_a\colon \ \mathbb{P} G' \to {\mathcal P} (a)\right\}_{a\in G_0}.
\end{gather*}
For any node $a'$ of $G'$ we have the canonical projection
\begin{gather*}
\pi'_{a'}\colon \ \mathbb{P} G'\to {\mathcal P}' (a').
\end{gather*}
We therefore def\/ine
\begin{gather*}
(\mathbb{P} \varphi)_a:= \pi'_{\varphi_a}\colon \ \mathbb{P} G' \to {\mathcal P}' (\varphi(a)) = {\mathcal P} (a)
\end{gather*}
for all $a\in G_0$.
By the universal property of the product $\mathbb{P} G= \bigsqcap\limits_{a\in G_0} {\mathcal P}(a)$ this def\/ines the desired
map $\mathbb{P} \varphi\colon \mathbb{P} G'\to \mathbb{P} G$.

The universal property of products also implies that the map $\mathbb{P}$ on morphisms of $\mathsf{Graph}/\mathsf{Man}$
as def\/ined above is actually a~functor.
That is,
\begin{gather*}
\mathbb{P} (\psi \circ \varphi) = \mathbb{P} \varphi\circ \mathbb{P} \psi
\end{gather*}
for any pair $(\psi, \varphi)$ of composable morphisms in $\mathsf{Graph}/\mathsf{Man}$.
\end{proof}

\begin{Remark}
Proposition~\ref{prop2.13} is an instance of a~category-theoretic result that holds in greater generality.
Namely, given a~category $\mathscr{C}$ with f\/inite products consider the category $\mathsf{FinSet}/\mathscr{C}$ whose
objects are pairs $(X, P)$, where~$X$ is a~f\/inite set and~$P$ is a~function that assignes to each element of~$X$ an object of $\mathscr{C}$.
The morphisms are commuting triangles.
There is a~contravariant functor $\mathbb{P}: (\mathsf{FinSet}/\mathscr{C})^{{\mathsf{op}}} \to \mathscr{C}$ which on objects is
given~by
\begin{gather*}
\mathbb{P}(X,P) = \bigsqcap_{x\in X} P(x).
\end{gather*}
\end{Remark}

\begin{Example}
\label{ex:1.0}
Suppose~$G$ is a~graph with two nodes $a$, $b$ and no edges, $G'$ is a~graph with one node $\{c\}$ and no edges,
${\mathcal P}'(c)$ is a~manifold~$M$, and $\varphi\colon G\to G'$ is the only possible map of graphs (it sends both
nodes to~$c$).
Suppose further that ${\mathcal P}\colon G_0\to \mathsf{Man}$ is given by ${\mathcal P} (a) = M = {\mathcal P}(b)$ (so
that ${\mathcal P}'\circ \varphi = {\mathcal P}$).
Then $\mathbb{P} G'\simeq M$,
\begin{gather*}
\mathbb{P} G = \{(x_a, x_b) \,|\, x_a\in {\mathcal P} (a), x_b \in {\mathcal P} (b)\} \simeq M\times M,
\end{gather*}
and $\mathbb{P}\varphi\colon M\to M\times M$ is the unique map with $(\mathbb{P} \varphi (x))_a = x$ and $(\mathbb{P}
\varphi (x))_b = x$ for all $x\in \mathbb{P} G'$.
Thus $\mathbb{P} \varphi\colon M\to M\times M$ is the diagonal map $x\mapsto (x,x)$.
\end{Example}

\begin{Example}
Let $(G,{\mathcal P})$, $(G', {\mathcal P}')$ be as in Example~\ref{ex:1.0} above and $\psi\colon (G',{\mathcal P}')\to
(G,{\mathcal P})$ be the map that sends the node~$c$ to~$a$.
Then $\mathbb{P} \psi\colon \mathbb{P} G \to \mathbb{P} G'$ is the map that sends $(x_a, x_b)$ to $x_a$.
\end{Example}

\begin{Remark}
\label{rmrk:2.16}
If $(G,{\mathcal P})$ is a~graph with a~phase function, that is, an object of $\mathsf{Graph}/\mathsf{Man}$, and
$\varphi\colon H\to G$ a~map of graphs then ${\mathcal P}\circ \varphi\colon H\to \mathsf{Man}$ is a~phase function and
$\varphi\colon (H,{\mathcal P}\circ \varphi)\to (G,{\mathcal P})$ is a~morphism in $\mathsf{Graph}/\mathsf{Man}$.
We then have a~map of manifolds
\begin{gather*}
\mathbb{P} \varphi\colon \ \mathbb{P} (H, {\mathcal P}\circ \varphi)\to \mathbb{P} (G,{\mathcal P}).
\end{gather*}
Similarly, a~commutative diagram
\begin{gather*}
\xy
(-8, 4)*+{ K} ="1";
(8, 4)*+{ H} ="2";
(0, -4)*+{G}="3";
{\ar@{->}^{ j} "1";"2"};
{\ar@{->}_{ \psi} "1";"3"};
{\ar@{->}^{\varphi} "2";"3"};
\endxy
\end{gather*}
of maps of graphs and a~phase space function ${\mathcal P}\colon G\to \mathsf{Man}$ give rise to the commutative diagram
of maps of manifolds
\begin{gather*}
\xy
(-18, 14)*+{\mathbb{P} (K, {\mathcal P}\circ \psi) } ="1";
(18, 14)*+{ \mathbb{P} (H, {\mathcal P}\circ \varphi)} ="2";
(0, -4)*+{\mathbb{P}(G,{\mathcal P})}="3";
{\ar@{->}_{ \mathbb{P} j} "2";"1"};
{\ar@{->}^{ \mathbb{P} \psi} "3";"1"};
{\ar@{->}_{\mathbb{P} \varphi} "3";"2"};
\endxy
\end{gather*}

\end{Remark}

\subsection{Embeddings and submersions from maps of graphs}

As we said in the introduction, the main goal of this paper is to construct maps of dynamical systems from graph
f\/ibrations.
In Proposition~\ref{prop2.13} we showed that a~map of networks $\varphi\colon (G,{\mathcal P})\to (G', {\mathcal P}')$
def\/ines a~map of manifolds $\mathbb{P} \varphi \colon \mathbb{P} (G', {\mathcal P}') \to \mathbb{P} (G,{\mathcal P})$.
In this subsection we prove that:
\begin{enumerate}\itemsep=0pt
\item[1)] If the map of graphs $\varphi\colon G\to G'$ is injective on nodes, then $\mathbb{P} \varphi$ is a~surjective
submersion,
\item[2)] if the map of graphs $\varphi\colon G\to G'$ is surjective on nodes, then $\mathbb{P} \varphi$ is an embedding.
\end{enumerate}
(Recall that a~smooth map between two manifolds is a~submersion if its dif\/ferential is onto at every point.
A~smooth map between two manifolds is an embedding if it is 1-1, its dif\/ferential is 1-1 everywhere and it is
a~homeomorphism onto its image.) Combined with Theorem~\ref{thm:main} below, this shows that surjective f\/ibrations of
networks of manifolds give rise to invariant dynamical subsystems and injective f\/ibrations give rise to projections of
dynamical systems.

\begin{Lemma}
\label{lemma:2.surj}
Suppose $\varphi\colon (G,{\mathcal P})\to (G',{\mathcal P}')$ is a~map of networks of manifolds such that the map on
nodes, $\varphi_0\colon G_0 \to G'_0$, is surjective.
Then $\mathbb{P}\varphi\colon \mathbb{P} G'\to \mathbb{P} G$ is an embedding whose image is the ``polydiagonal''
\begin{gather*}
\Delta_\varphi = \{x\in \mathbb{P} G\,|\, x_a = x_b\ \text{whenever}\ \varphi(a) = \varphi(b)\}.
\end{gather*}
\end{Lemma}
\begin{proof}
Assume f\/irst for simplicity that $G'$ has only one vertex $*$ and ${\mathcal P}'(*) = M$.
Then for any vertex~$a$ of~$G$ we have
\begin{gather*}
{\mathcal P}(a) = {\mathcal P}' (\varphi(a)) = {\mathcal P}'(*) = M,
\end{gather*}
$\mathbb{P} G' = M$ and $\mathbb{P} G = M\times \dots \times M$ ($|G_0|$ copies), where as before $G_0$ is the set of vertices of the graph~$G$.
In this case the proof of Proposition~\ref{prop2.13} shows that the map $\mathbb{P}\varphi\colon M\to M^{G_0}$ is of the form
\begin{gather*}
\mathbb{P} \varphi (x) = (x,\ldots, x)
\end{gather*}
for all $x\in M$.
This is clearly an embedding.
In general,
\begin{gather*}
\mathbb{P} \varphi\colon \ \mathbb{P} G' =\bigsqcap_{a'\in G'_0} {\mathcal P}' (a') \to \bigsqcap_{a'\in G'_0}
\left(\bigsqcap_{a\in \varphi^{-1} (a')} {\mathcal P}(a)\right) = \mathbb{P} G
\end{gather*}
is the product of maps of the form
\begin{gather*}
{\mathcal P}' (a') \to \bigsqcap_{a\in \varphi^{-1} (a')} {\mathcal P}(a),
\qquad
x\mapsto (x,\ldots, x).\tag*{\qed}
\end{gather*}
\renewcommand{\qed}{}
\end{proof}

\begin{Lemma}
\label{lemma:2.inj}
Suppose $\varphi\colon (G,{\mathcal P})\to (G',{\mathcal P}')$ is a~map of networks of manifolds such that the map
$\varphi_0\colon G_0 \to G'_0$ on nodes is injective.
Then $\mathbb{P}\varphi\colon \mathbb{P} G'\to \mathbb{P} G$ is a~surjective submersion.
\end{Lemma}
\begin{proof}
Since $\varphi\colon G\to G'$ is injective, the set of nodes $G_0'$ of $G'$ can be partitioned as the disjoint union of
the image $\varphi(G_0)$, which is a~copy of $G_0$, and the complement.
Hence
\begin{gather*}
\mathbb{P} G' \simeq \bigsqcap_{a\in G_0} {\mathcal P}(\varphi(a)) \times \bigsqcap_{a'\not\in \varphi(G_0)} {\mathcal
P}' (a') \simeq \mathbb{P} G \times \bigsqcap_{a'\not\in \varphi(G_0)} {\mathcal P}' (a').
\end{gather*}
With respect to this identif\/ication of $\mathbb{P} G'$ with $\mathbb{P} G \times \bigsqcap\limits_{a'\not\in \varphi(G_0)}
{\mathcal P}' (a')$ the map $\mathbb{P}\varphi\colon \mathbb{P} G'\to \mathbb{P} G$ is the projection
\begin{gather*}
\mathbb{P} G \times \bigsqcap_{a'\not\in \varphi(G_0)} {\mathcal P}' (a') \to \mathbb{P} G.
\end{gather*}
which is a~surjective submersion.
\end{proof}

\subsection{Open systems and their interconnections}

\noindent Having set up a~consistent way of assigning phase spaces to graphs, we now take up continuous time dynamical
systems.
We start by recalling a~def\/inition of an open (control) systems, which is essentially due to Brockett~\cite{Brockett}.
It is not the most general def\/inition; it is more than enough for this paper.

\begin{Definition}
\label{def:open-sys}
A~{\it continuous time control system} (or an {\it open system}) on a~manifold~$M$ is a~surjective submersion $p\colon
Q\to M$ from some manifold~$Q$ together with a~smooth map $ F\colon Q \to TM$ so that
\begin{gather*}
F(q) \in T_{p(q)} M
\end{gather*}
for all $q\in Q $.
That is, the diagram $
\xy (-10, 6)*+{Q}="1";
(6, 6)*+{TM} ="2";
(6,-3)*+{M}="3";
{\ar@{->}_{ p}  "1";"3"};
{\ar@{->}^{F} "1";"2"};
{\ar@{->}^{\pi} "2";"3"}; \endxy
$ commutes.
Here $\pi\colon TM \to M$ is the canonical projection.
\end{Definition}

\begin{Definition}[$\mathsf{Control} (M\times U\to M)$] Given a~manifold~$M$ of ``state variables'' and a~manifold~$U$ of ``control
variables'' we may consider control systems of the form
\begin{gather*}
F\colon  \ M\times U\to TM,
\qquad
F(x,u)\in T_x M
\qquad
\text{for all}
\quad
(x,u)\in M\times U.
\end{gather*}
The collection of all such control systems forms a~vector space $\mathsf{Control} (M\times U\to M)$.
Explicitly
\begin{gather*}
\mathsf{Control} (M\times U\to M):= \{F\colon M\times U\to TM \,|\, F(x,u)\in T_x M \ \text{for all} \ (x,u)\in M\times U\}.
\end{gather*}
\end{Definition}

Now suppose we are given a~f\/inite family $\{F_i\colon M_i \times U_i \to TM_i\}_{i=1}^N$ of control systems and we want
to somehow interconnect them to obtain a~closed system $\mathscr{I} (F_1, \ldots, F_N)$, that is, a~vector f\/ield on the
product $\bigsqcap_i M_i$.
What additional data do we need to def\/ine the interconnection map
\begin{gather*}
\mathscr{I}\colon \ \bigsqcap_i \mathsf{Control}(M_i\times U_i \to M_i) \to \Gamma \bigg(T\bigg(\bigsqcap_i M_i\bigg)\bigg)?
\end{gather*}
An answer is given by the following proposition:

\begin{Proposition}
\label{prop:3.2}
Given a~family $\{p_j \colon M_j \times U_j \to M_j\}_{j=1}^N$ of projections on the first factor and a~family of smooth
maps $\{s_j\colon \bigsqcap M_i \to M_j\times U_j\}$ so that the diagrams
\begin{gather*}
\xy (-10, 6)*+{M_j \times U_j}="1";
(-10,-6)*+{\bigsqcap M_i}="3";
(10,-6)*+{ M_j}="4";
{\ar@{->}^{ p_j}  "1";"4"};
{\ar@{->}^{s_j} "3";"1"};
{\ar@{->}_{pr_j} "3";"4"};
\endxy
\end{gather*}
commute for each index~$j$, there is an interconnection map $\mathscr{I}$ making the diagrams
\begin{gather*}
\xy
(-30, 6)*+{\bigsqcap_i \mathsf{Control }(M_i\times U_i \to M_i)}="1";
(30, 6)*+{\Gamma
(T(\bigsqcap_i M_i))} ="2";
(-30,-6)*+{\mathsf{Control }(M_j\times U_j \to M_j)}="3";
(30,-6)*+{\mathsf{Control } (\bigsqcap_i M_i \xrightarrow {pr_j} M_j)}="4";
{\ar@{->}^{\mathscr{I}} "1";"2"};
{\ar@{->}^{\varpi_j =D(pr_j)\circ -} "2";"4"};
{\ar@{->}^{} "1";"3"};
{\ar@{->}_{\mathscr{I}_j} "3";"4"};
\endxy
\end{gather*}
commute for each~$j$.
The components $\mathscr{I}_j$ of the interconnection map $\mathscr{I}$ are defined by $\mathscr{I}_j (F_j):= F_j \circ
s_j$ for all~$j$, where $D (pr_j)\colon T \bigsqcap M_i \to TM_j$ denotes the differential of the canonical projection
$pr_j\colon \bigsqcap M_i \to M_j$.
\end{Proposition}

\begin{proof}
The space of vector f\/ields $\Gamma \left(T\left(\bigsqcap_i M_i\right)\right)$ on the product $\bigsqcap_i M_i$ is the product of vector
spaces $\mathsf{Control} (\bigsqcap_i M_i \to M_j)$:
\begin{gather*}
\Gamma \bigg(T\bigg(\bigsqcap_i M_i\bigg)\bigg) = \bigsqcap_j \mathsf{Control} \bigg(\bigsqcap_i M_i \xrightarrow{pr_j}M_j\bigg).
\end{gather*}
In other words a~vector f\/ield $X $ on the product $\bigsqcap_i M_i$ is a~tuple $X = (X_1, \ldots, X_N)$, where
\begin{gather*}
X_j:= D (pr_j)\circ X.
\end{gather*}
Each component $X_j \colon \bigsqcap_i M_i \to TM_i$ is a~control system.

To def\/ine a~map from a~vector space into a~product of vector spaces it is enough to def\/ine a~map into each of the factors.
We have canonical projections
\begin{gather*}
\pi_j\colon \ \bigsqcap_i \mathsf{Control}(M_i\times U_i \to M_i) \to \mathsf{Control}(M_j\times U_j\to M_j),
\qquad
j=1,\ldots, N.
\end{gather*}
Consequently to def\/ine the interconnection map $\mathscr{I}$ it is enough to def\/ine the maps
\begin{gather*}
\mathscr{I}_j\colon \ \mathsf{Control}(M_j\times U_j\to M_j) \to \mathsf{Control} \bigg(\bigsqcap_i M_i \xrightarrow {pr_j}M_j\bigg).
\end{gather*}
for each index~$j$.
We therefore def\/ine the maps $\mathscr{I}_j\colon \mathsf{Control}(M_j\times U_j\to M_j) \to \mathsf{Control}
(\bigsqcap_i M_i \xrightarrow {pr_j} M_j)$, $1\leq j\leq N$,~by
\begin{gather*}
\mathscr{I}_j (F_j):= F_j \circ s_j. \tag*{\qed}
\end{gather*}
\renewcommand{\qed}{}
\end{proof}

\begin{Remark}
\label{rmrk:2.21}
It will be useful for us to remember that the canonical projections
\begin{gather*}
\varpi_j\colon \ \Gamma \Big(T\bigsqcap M_i\Big)
\to \mathsf{Control} \Big(\bigsqcap M_i \to M_j\Big)
\end{gather*}
are given~by
\begin{gather*}
\varpi_j (X) = D (pr_j)\circ X,
\end{gather*}
where as before $D (pr_j)\colon T \bigsqcap M_i \to TM_j$ are the dif\/ferentials of the canonical projections $pr_j\colon\bigsqcap M_i \to M_j$.
\end{Remark}

\subsection{Interconnections and graphs}

We next explain how f\/inite directed graphs whose nodes are decorated with phase spaces, that is, networks of manifolds
in the sense of Def\/inition~\ref{def.network.mfld} give rise to interconnection maps.
To do this precisely it is useful to have a~notion of {\em input trees} of a~directed graph.
This notion is a~generalization of the notion of an {\em input set} of Golubitsky~et al.\ ({\em op.\ cit.})~\cite{Golubitsky.Stewart.Torok.05,Stewart.Golubitsky.Pivato.03}.
Given a~graph, an input tree~$I(a)$ of a~vertex~$a$ is~-- roughly~-- the vertex itself and all of the arrows leading into it.
We want to think of this as a~graph in its own right, as follows.

\begin{Definition}[input tree]\label{def:input-tree}
Given a~vertex~$a$ of a~graph $ G$ we def\/ine the {\em input tree} $I(a)$ to be a~graph with the set of vertices $I(a)_0$
given~by
\begin{gather*}
I(a)_0:= \{a\}\sqcup \mathfrak t^{-1} (a),
\end{gather*}
where, as before, the set $\mathfrak t^{-1} (a)$ is the set of arrows in $ G$ with target~$a$.
The set of edges $I(a)_1$ of the input tree is the set of pairs
\begin{gather*}
I(a)_1:= \{(a, \gamma)\,|\, \gamma \in G_1, \mathfrak t(\gamma) =a \},
\end{gather*}
and the source and target maps $I(a)_1\rightrightarrows I(a)_0$ are def\/ined~by
\begin{gather*}
\mathfrak s(a, \gamma) = \gamma
\qquad
\text{and}
\qquad
\mathfrak t(a, \gamma) = a.
\end{gather*}
In pictures,

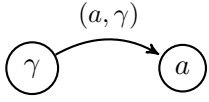
\begin{tikzpicture}
[->,>=stealth',shorten >=1pt,auto,node distance=2cm, thick,main node/.style={circle,fill=white,draw}]

\node[main node] at (2,0) (1) {$\gamma$}; \node[main node] at (4,0) (2) {$a$};

\path[every node/.style={font=\sffamily\small}] (1) edge [bend left] node {$(a,\gamma)$} (2);
\end{tikzpicture}
\end{Definition}

\begin{Example}
Consider the graph

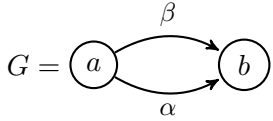
\begin{tikzpicture}
[->,>=stealth',shorten >=1pt,auto,node distance=2cm, thick,main node/.style={circle,draw}]

\node at (0,-1) {$G = $}; \node[main node] at (0.8,-1) (1) {$a$}; \node[main node] at (2.8,-1) (2) {$b$};

\path[every node/.style={font=\sffamily\small}] (1) edge [bend left] node {$\beta$} (2) (1) edge [bend right] node
[below] {$\alpha$} (2);
\end{tikzpicture}

\noindent
as in Example~\ref{example:1}.
Then the input tree $I(a)$ is the graph with one node~$a$ and no edges: $I(a) = \{a\}$ (see
Notation~\ref{notation:2.11}).
The input tree $I(b)$ has three nodes and two edges:

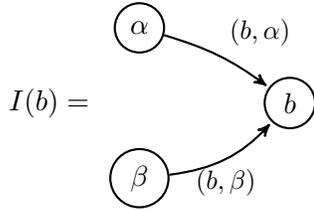
\begin{tikzpicture}
[->,>=stealth',shorten >=1pt,auto,node distance=2cm, thick,main
node/.style={circle,fill=white,draw,font=\sffamily\bfseries}]

\node at (0,0) {$I(b) = $}; \node[main node] at (1.2,1) (1) {$\alpha$}; \node[main node] at (1.2,-1) (2) {$\beta$};
\node[main node] at (3.2,0) (3) {$b$};

\path[every node/.style={font=\sffamily\small}] (1) edge [bend left = 10] node {$(b,\alpha)$} (3) (2) edge [bend right =
20] node [below] {$(b,\beta)$} (3);
\end{tikzpicture}

Notice that our def\/inition of input tree ``pulls apart'' multiple edges coming from a~common vertex.
\end{Example}

\begin{Remark}
\label{remark:xi}
For each node~$a$ of a~graph~$G$ we have a~natural map of graphs
\begin{gather*}
\xi= \xi_a \colon  \ I(a)\to G.
\end{gather*}
It is def\/ined by sending the edge of the form $\gamma \xrightarrow{(a,\gamma)} a$ to the edge $\mathfrak s (\gamma)\xrightarrow{\gamma} a$.
Note that the map $\xi$ need not be injective on vertices.
\end{Remark}

\begin{Proposition}
\label{prop:2.26}
Given a~graph~$G$ with a~phase space function ${\mathcal P}\colon G_0 \to \mathsf{Man}$, that is, a~network
$(G,{\mathcal P})$ of manifolds, we have commutative diagrams of maps of manifolds
\begin{gather*}
\xy
(-40, 10)*+{ \mathbb{P} I(a)} ="0";
(-10, 10)*+{
{\mathcal P} (a)\times \bigsqcap\limits_{\gamma \in \mathfrak t^{-1} (a)}{\mathcal P} (\mathfrak s (\gamma))} ="1";
(40, 10)*+{ \mathbb{P} a} ="4";
(30, 10)*+{ {\mathcal P}(a)} ="2";
(-10, -10)*+{\bigsqcap\limits_{b\in G_0} {\mathcal P}(b)}="3";
(-30, -10)*+{\mathbb{P} G=};
{\ar@{=}^{ } "1";"0"};{\ar@{=}^{ } "2";"4"};
{\ar@{->}^(.75){ \mathbb{P} j_a} "1";"2"};
{\ar@{->}_{ \mathbb{P} \iota_a} "3";"2"};
{\ar@{->}^{\mathbb{P} \xi} "3";"1"};
\endxy
\end{gather*}
for each node~$a$ of the graph~$G$.
\end{Proposition}

\begin{proof}
Let~$a$ be a~node of the graph~$G$.
We then have a~graph $\{a\}$ with one node and no arrows.
Denote the inclusion of $\{a\}$ in~$G$ by $\iota_a$ and the inclusion into its input tree $I(a)$ by $j_a$.
Then the diagram of maps of graphs
\begin{gather*}
\xy
(-8, 10)*+{ \{a\}} ="1";
(8, 10)*+{ I(a)} ="2";
(0, -2)*+{G}="3";
{\ar@{->}^{ j_a} "1";"2"};
{\ar@{->}_{ \iota_a} "1";"3"};
{\ar@{->}^{\xi} "2";"3"};
\endxy
\end{gather*}
commutes.
By Remark~\ref{rmrk:2.16} we have a~commuting diagram of maps of manifolds
\begin{gather*}
\xy
(-10, 10)*+{ \mathbb{P}\{a\}} ="1";
(10, 10)*+{ \mathbb{P} I(a)} ="2";
(0, -5)*+{\mathbb{P} G}="3";
{\ar@{<-}^{ \mathbb{P} j_a} "1";"2"};
{\ar@{<-}_{ \mathbb{P} \iota_a} "1";"3"};
{\ar@{<-}^{\mathbb{P} \xi} "2";"3"};
\endxy
\end{gather*}

Let us now examine more closely the map $\mathbb{P} j_a\colon \mathbb{P} I(a)\to \mathbb{P} a$.

Since the set of nodes $I(a)_0$ of the input tree $I(a)$ is the disjoint union
\begin{gather*}
I(a)_0 = \{a\} \sqcup \mathfrak t^{-1} (a),
\end{gather*}
and since $\xi_a (\gamma) = \mathfrak s (\gamma)$ for any $\gamma \in \mathfrak t^{-1} (a) \subset I(a)_0$, we have
\begin{gather*}
\mathbb{P} I(a) = {\mathcal P}(a) \times \bigsqcap_{\gamma \in \mathfrak t^{-1} (a)} {\mathcal P} (\mathfrak s(\gamma)).
\end{gather*}
Since $j_a\colon \{a\} \to I(a)_0 = \{a\} \sqcup \mathfrak t^{-1} (a)$ is the inclusion,
\begin{gather*}
\mathbb{P} j_a \colon \ \mathbb{P} I(a) \to \mathbb{P} a
\end{gather*}
is the projection
\begin{gather*}
{\mathcal P}(a) \times \bigsqcap_{\gamma \in \mathfrak t^{-1} (a)} {\mathcal P} (\mathfrak s(\gamma)) \to \mathbb{P} a.
\end{gather*}
Similarly
\begin{gather*}
\mathbb{P} \iota_a\colon \ \mathbb{P} G \to \mathbb{P} a
\end{gather*}
is the projection
\begin{gather*}
\bigsqcap_{b\in G_0} {\mathcal P} (b) \to {\mathcal P}(a).
\end{gather*}
The proposition follows from these two observations.
\end{proof}

\begin{Example}\label{example:3}
Suppose
\begin{gather*}
\begin{tikzpicture}
[->,>=stealth',shorten >=1pt,auto,node distance=2cm, thick,main node/.style={circle,draw,font=\sffamily\bfseries}]
\node at (1,-1) {$G = $}; \node[main node] at (1.9,-1) (1) {a}; \node[main node] at (3.9,-1) (2) {b};
\path[every node/.style={font=\sffamily\small}] (1) edge [bend left] node {} (2) (1) edge [bend right] node [below] {}
(2);
\end{tikzpicture}
\end{gather*}
is a~graph as in Example~\ref{example:1} and suppose ${\mathcal P}\colon G_0 \to \mathsf{Man}$ is a~phase space function.
Then
\begin{gather*}
\mathbb{P} I(b) \simeq {\mathcal P} (a) \times {\mathcal P} (a) \times {\mathcal P} (b),
\end{gather*}
$\mathbb{P} j_b$ is the projection ${\mathcal P} (a) \times {\mathcal P} (a) \times {\mathcal P} (b) \to {\mathcal P}
(b)$, and
\begin{gather*}
\mathsf{Control} (\mathbb{P} I(b) \to \mathbb{P} b) = \mathsf{Control} ({\mathcal P} (a) \times {\mathcal P} (a) \times
{\mathcal P} (b) \to {\mathcal P} (b)).
\end{gather*}
On the other hand $\mathbb{P} I(a) = {\mathcal P} (a)$, $\mathbb{P} j_a\colon {\mathcal P}(a)\to {\mathcal P} (a)$ is
the identity map and
\begin{gather*}
\mathsf{Control} (\mathbb{P} I(a)\to \mathbb{P} a) = \Gamma (T{\mathcal P}(a)),
\end{gather*}
the space of vector f\/ields on the manifold ${\mathcal P} (a)$.
\end{Example}

\begin{Notation}
Given a~network $(G,{\mathcal P})$ of manifolds we have a~product of vector spaces
\begin{gather*}
\mathpzc{Ctrl}(G,{\mathcal P}):= \bigsqcap_{a\in G_0} \mathsf{Control}(\mathbb{P} I(a)\to \mathbb{P} a).
\end{gather*}
The elements of $\mathpzc{Ctrl}(G,{\mathcal P})$ are unordered tuples of $(w_a)_{a\in G_0}$ of control systems (q.v.\
Lem\-ma~\ref{empt:cat-prod-vect}).
We may think of them as sections of the vector bundle $\bigsqcup\limits_{a\in G_0}\mathsf{Control}(\mathbb{P} I(a)\to
\mathbb{P} a) \to G_0$ over the vertices of~$G$.
\end{Notation}
It is easy to see that Propositions~\ref{prop:3.2} and~\ref{prop:2.26} give us
\begin{Theorem}
\label{thm:3.8}
Given a~network $(G,{\mathcal P})$ of manifolds, there exists a~natural interconnection map
\begin{gather*}
\mathscr{I}\colon \ \bigsqcap_{a\in G_0} \mathsf{Control} (\mathbb{P} I(a)\to \mathbb{P} a) \to \Gamma (T\mathbb{P} G)
\end{gather*}
with
\begin{gather*}
\varpi_a \circ \mathscr{I} ((w_b)_{b\in G_0}) = w_a \circ \mathbb{P} j_a
\end{gather*}
for all nodes $a\in G_0$.
Here $\varpi_a\colon \Gamma (T\mathbb{P} G) \to \mathsf{Control}(\mathbb{P} G_0 \xrightarrow{\mathbb{P} \iota_a}
\mathbb{P} a)$ are the projection maps; $\varpi_a = D (\mathbb{P} \iota_a)$ $($q.v.~Remark~{\rm \ref{rmrk:2.21})}.
\end{Theorem}

\begin{Example}
Consider the graph~$G$ as in Examples~\ref{example:1} and~\ref{example:3} with a~phase space function ${\mathcal
P}\colon G_0\to \mathsf{Man}$.
Then the vector f\/ield
\begin{gather*}
X = \mathscr{I} (w_a, w_b)\colon \ {\mathcal P} (a)\times {\mathcal P} (b) \to T{\mathcal P} (a)\times T{\mathcal P} (b)
\end{gather*}
is of the form
\begin{gather*}
X(x,y) = (w_a (x), w_b (x,x, y))
\qquad
\text{for all}
\quad
(x,y) \in {\mathcal P} (a)\times {\mathcal P} (b).
\end{gather*}
\end{Example}
\begin{Example}
Consider the graph
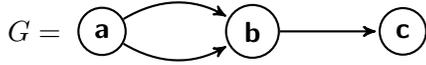
\begin{gather*}
\begin{tikzpicture}
[->,>=stealth',shorten >=1pt,auto,node distance=2cm, thick,main node/.style={circle,draw,font=\sffamily\bfseries}]
\node at (1,-1) {$G = $}; \node[main node] at (1.9,-1) (1) {a}; \node[main node] at (3.9,-1) (2) {b}; \node[main node] at
(5.9,-1) (3) {c};
\path[every node/.style={font=\sffamily\small}] (1) edge [bend left] node {} (2) (1) edge [bend right] node [below] {}
(2) (2) edge node {}(3);
\end{tikzpicture}
\end{gather*}
and let ${\mathcal P} \colon G_0\to \mathsf{Man}$ be a~phase space function.
Then
\begin{gather*}
\left(\mathscr{I} (w_a, w_b, w_c)\right) (x,y,z) = (w_a (x), w_b (x,x, y), w_c (y,z))
\end{gather*}
for all $(w_a,w_b, w_c) \in \mathpzc{Ctrl} (G,{\mathcal P})$ and all $(x,y,z)\in {\mathcal P} (a)\times {\mathcal P} (b)
\times {\mathcal P} (c)$.
\end{Example}

\section{Maps of dynamical systems from f\/ibrations}

Following Boldi and Vigna~\cite{Vigna1} (see also~\cite{VignaFibPage}) we single out a~class of maps of graphs called
graph f\/ibrations.
\begin{Definition}
\label{graph fibration}
A~map $\varphi\colon G\to G'$ of directed graphs is a~{\it graph fibration} if for any vertex~$a$ of~$G$ and any edge~$e'$ of $ G'$ ending at $\varphi(a)$ there is a~unique edge~$e$ of $ G$ ending at~$a$ with $\varphi (e) = e'$.
\end{Definition}

\begin{Example}
\label{example:5}
The map of graphs
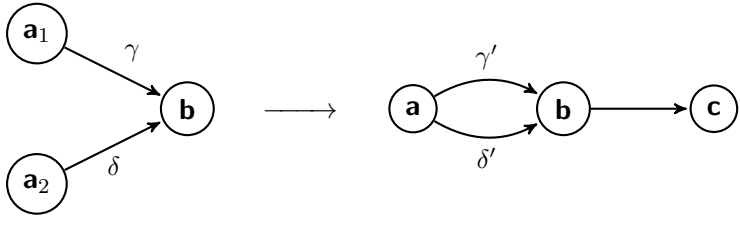
\begin{gather*}
\begin{tikzpicture}
[->,>=stealth',shorten >=1pt,auto,node distance=2cm, thick,main node/.style={circle,draw,font=\sffamily\bfseries}]
\node[main node] at (-2,0) (a1) {a$_1$}; \node[main node] at (-2,-2) (a2) {a$_2$}; \node[main node] at (0,-1) (bb) {b};
\node[main node] at (3,-1) (1) {a}; \node[main node] at (5,-1) (2) {b}; \node[main node] at (7,-1) (3) {c}; \node at
(1.5,-1) (arrow) {$\xrightarrow{
\qquad
}$};
\path[every node/.style={font=\sffamily\small}] (a1) edge node {$\gamma$} (bb) (a2) edge node [below] {$\delta$} (bb)
(1) edge [bend left] node {$\gamma'$} (2) (1) edge [bend right] node [below] {$\delta'$} (2) (2) edge node {}(3)
;
\end{tikzpicture}
\end{gather*}
sending the edge~$\gamma$ to $\gamma'$ and the edge~$\delta$ to $\delta'$ is a~graph f\/ibration.
\end{Example}

\begin{Remark}
\label{empt:4.2}
Given any map $\varphi\colon G\to G'$ of graphs and a~node~$a$ of~$G$ there is an induced map of input trees
\begin{gather*}
\varphi_a\colon \  I(a) \to I(\varphi(a)).
\end{gather*}
On edges of $I(a)$ the map is def\/ined~by
\begin{gather*}
\varphi (a, \gamma):= (\varphi (a), \varphi (\gamma))
\end{gather*}
(cf.\ Def\/inition~\ref{def:input-tree}).
Moreover the diagram of graphs
\begin{gather*}
\xy
(-12, 8)*+{I(a)} ="1";
(12, 8)*+{I(\varphi(a))} ="2";
(-12, -10)*+{G}="3";
(12, -10)*+{G'}="4";
{\ar@{->}^{ \varphi_a} "1";"2"};
{\ar@{->}_{ \xi_a} "1";"3"};
{\ar@{->}^{ \xi_{\varphi(a)}} "2";"4"};
{\ar@{->}^{\varphi} "3";"4"};
\endxy
\end{gather*}
commutes (the map $\xi_a\colon I(a)\to G$ from an input tree to the original graph is def\/ined in
Remark~\ref{remark:xi}).
\end{Remark}

\begin{Lemma}
\label{lemma:4.3}
If $\varphi\colon G\to G'$ is a~graph fibration then the induced maps
\begin{gather*}
\varphi_a\colon \ I(a)\to I(\varphi(a))
\end{gather*}
of input trees defined above are isomorphisms for all nodes~$a$ of~$G$.
\end{Lemma}

\begin{proof}
Given an edge $(\varphi(a), \gamma')$ of $I(\varphi(a))$ there is a~unique edge~$\gamma$ of~$G$ with $\varphi(\gamma) =
\gamma'$ and $\mathfrak t(\gamma) = a$ and consequently $\varphi_a (a,\gamma) = (\varphi(a),\gamma')$.
It follows that $\varphi_a$ is bijective on vertices and edges.
\end{proof}
\begin{Remark}
The converse is true as well: if the induced maps $\varphi_a\colon I(a)\to I(\varphi(a))$ are isomorphisms for all
nodes~$a$ of~$G$ then $\varphi\colon G\to G'$ is a~graph f\/ibration.
\end{Remark}

Recall that a~map from a~network $(G,{\mathcal P})$ to a~network $(G',{\mathcal P}')$ is a~map of graphs $\varphi\colon
G\to G'$ with the property that
\begin{gather*}
{\mathcal P}'\circ \varphi = {\mathcal P}.
\end{gather*}

\begin{Definition}[f\/ibration of networks of manifolds]
A~map of networks $\varphi\colon (G,{\mathcal P}) \to (G',{\mathcal P}')$ of
manifolds is a~{\it fibration} if $\varphi\colon G\to G'$ is a~graph f\/ibration.
\end{Definition}
\begin{Remark}[the category $(\mathsf{Man}/\mathsf{Graph})_{\mbox{\sf{\tiny{f\/ib}}}} $ of networks of manifolds and f\/ibrations] We note
that the composit of two f\/ibrations is again a~f\/ibration.
Consequently networks of manifolds and f\/ibrations form a~category which we denote~by
$(\mathsf{Man}/\mathsf{Graph})_{\mbox{\sf{\tiny{f\/ib}}}} $.
\end{Remark}
Theorem~\ref{thm:4.4} below is our reason for singling out f\/ibrations of networks.
\begin{Theorem}
\label{thm:4.4}
A~fibration $\varphi\colon (G,{\mathcal P})\to (G',{\mathcal P}')$ of networks induces a~linear map
\begin{gather*}
\varphi^*\colon \ \mathpzc{Ctrl}(G',{\mathcal P}')\to \mathpzc{Ctrl}(G,{\mathcal P}).
\end{gather*}
\end{Theorem}

\begin{proof}
Since
\begin{gather*}
\mathpzc{Ctrl}(G,{\mathcal P}) = \bigsqcap_{a\in G_0} \mathsf{Control} (\mathbb{P} I(a)\to \mathbb{P} a)
\end{gather*}
is a~product of vector spaces, the map $\varphi^*$ is uniquely determined by maps from $\mathpzc{Ctrl}(G',{\mathcal
P}')$ to the factors $\mathsf{Control} (\mathbb{P} I(a)\to \mathbb{P} a)$, $a\in G_0$.
On the other hand we have canonical projections
\begin{gather*}
\pi_b\colon \ \mathpzc{Ctrl}(G',{\mathcal P}') = \bigsqcap_{c\in G'_0} \mathsf{Control} (\mathbb{P} I(c)\to \mathbb{P}
c)\to \mathsf{Control} (\mathbb{P} I(b)\to \mathbb{P} b)
\end{gather*}
for all $b\in G_0'$.
Hence in order to def\/ine the map $\varphi^*$ it is enough to def\/ine maps of vector spaces
\begin{gather*}
\varphi_a^*\colon \ \mathsf{Control} (\mathbb{P} I(\varphi(a))\to \mathbb{P} \varphi(a))\to \mathsf{Control} (\mathbb{P}
I(a)\to \mathbb{P} a)
\end{gather*}
for all nodes~$a$ of the graph~$G$.
By Remark~\ref{empt:4.2} the diagram
\begin{gather*}
\xy
(-12, 8)*+{I(a)} ="1";
(12, 8)*+{I(\varphi(a))} ="2";
(-12, -10)*+{G}="3";
(12, -10)*+{G'}="4";
(0, -22)*+{\mathsf{Man}}="5";
{\ar@{->}^{ \varphi_a} "1";"2"};
{\ar@{->}_{ \xi_a} "1";"3"};
{\ar@{->}^{ \xi_{\varphi(a)}} "2";"4"};
{\ar@{->}^{\varphi} "3";"4"};
{\ar@{->}_{{\mathcal P}} "3";"5"};
{\ar@{->}^{{\mathcal P}'} "4";"5"};
\endxy
\end{gather*}
commutes for each $a\in G_0$.
Let
\begin{gather*}
\varphi|_{\{a\}}\colon \ \{a\} \to \{\varphi (a)\}
\end{gather*}
denote the restriction of $\varphi\colon G\to G'$ to the subgraph $\{a\}\hookrightarrow G$.
It is easy to see that the diagrams
\begin{gather*}
\xy
(-12, 8)*+{I(a)} ="1";
(12, 8)*+{I(\varphi(a))} ="2";
(-12, -10)*+{\{a\}}="3";
(12, -10)*+{\{\varphi(a)\}}="4";
{\ar@{->}^{ \varphi_a} "1";"2"};
{\ar@{->}_{ j_a} "3";"1"};
{\ar@{->}^{ j_{\varphi(a)}} "4";"2"};
{\ar@{->}^{\varphi|_{\{a\}}} "3";"4"};
\endxy
\end{gather*}
commutes as well.
By Lemma~\ref{lemma:4.3} the map $\varphi_a$ is an isomorphism of graphs.
Hence
\begin{gather*}
\mathbb{P} \varphi_a\colon  \ \mathbb{P} I(a) \to \mathbb{P} I(\varphi(a))
\end{gather*}
is an isomorphism of manifolds.
Def\/ine
\begin{gather*}
\varphi_a^*\colon \ \mathsf{Control} (\mathbb{P} I(\varphi(a))\to \mathbb{P} \varphi(a))\to \mathsf{Control} (\mathbb{P}
I(a)\to \mathbb{P} a)
\end{gather*}
by
\begin{gather*}
\varphi_a^*(F) = D\mathbb{P} (\varphi|_{\{a\}}) \circ F \circ (\mathbb{P} \varphi_a)^{-1}
\end{gather*}
for all $F\in \mathsf{Control} (\mathbb{P} I(\varphi(a)))$.
By the universal property of products this gives us the desired map $\varphi^*$.
Moreover the diagrams
\begin{gather*}
\xy
(-30, 8)*+{\mathpzc{Ctrl}(G',{\mathcal P}')} ="1";
(30, 8)*+{\mathpzc{Ctrl}(G,{\mathcal P})} ="2";
(-30, -10)*+{\mathsf{Control }(\mathbb{P} I(\varphi(a))\to \mathbb{P}\varphi(a))}="3";
(30, -10)*+{\mathsf{Control }(\mathbb{P} I(a) \to \mathbb{P} a)
}="4";
{\ar@{->}^{ \varphi^*} "1";"2"};
{\ar@{->}_{ \pi_{\varphi(a)}} "1";"3"};
{\ar@{->}^{ \pi_a} "2";"4"};
{\ar@{->}^{\varphi_a^*} "3";"4"};
\endxy
\end{gather*}
commute for all $a\in G_0$.
\end{proof}

\begin{Example}\label{example:6}
We write down an example of the map $\varphi^*$ constructed in Theorem~\ref{thm:4.4}.
Consider the graph f\/ibration $\varphi\colon G\to G'$:
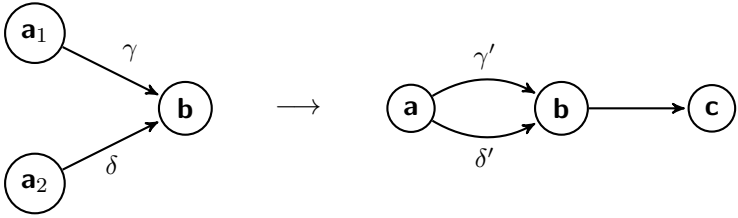
\begin{gather*}
\begin{tikzpicture}
[->,>=stealth',shorten >=1pt,auto,node distance=2cm, thick,main node/.style={circle,draw,font=\sffamily\bfseries}]
\node[main node] at (-2,0) (a1) {a$_1$}; \node[main node] at (-2,-2) (a2) {a$_2$}; \node[main node] at (0,-1) (bb) {b};
\node[main node] at (3,-1) (1) {a}; \node[main node] at (5,-1) (2) {b}; \node[main node] at (7,-1) (3) {c}; \node at
(1.5,-1) (arrow) {$\longrightarrow$};
\path[every node/.style={font=\sffamily\small}] (a1) edge node {$\gamma$} (bb) (a2) edge node [below] {$\delta$} (bb)
(1) edge [bend left] node {$\gamma'$} (2) (1) edge [bend right] node [below] {$\delta'$} (2) (2) edge node {}(3)
;
\end{tikzpicture}
\end{gather*}
as in Example~\ref{example:5}.
Let ${\mathcal P}'\colon G_0'\to \mathsf{Man}$ be a~phase space function.
Then
\begin{gather*}
\mathpzc{Ctrl}(G',{\mathcal P}') = \{(w_a\colon {\mathcal P} (a)\to T{\mathcal P} (a),w_b\colon {\mathcal P}'(a)\times
{\mathcal P}'(a) \times {\mathcal P}'(b)\to T{\mathcal P} (b),\\
\hphantom{\mathpzc{Ctrl}(G',{\mathcal P}') = \{}{}
w_c\colon {\mathcal P}(b)\times {\mathcal P} (c) \to T{\mathcal P}(c))\},
\\
\mathpzc{Ctrl}(G,{\mathcal P}'\circ \varphi) = \{(w'_{a_1}\colon {\mathcal P}' (a)\to {\mathcal P}' (a), w_{a_2}\colon
{\mathcal P}' (a)\to {\mathcal P}' (a),\\
\hphantom{\mathpzc{Ctrl}(G,{\mathcal P}'\circ \varphi) = \{}{}
w_b\colon {\mathcal P}'(a)\times {\mathcal P}'(a) \times {\mathcal P}'(b)\to T{\mathcal P}' (b))\}
\end{gather*}
and
\begin{gather*}
\varphi^* (w_a',w_b',w_c') = (w_a', w_a', w_b').
\end{gather*}
\end{Example}

\begin{Remark}[the category $\mathsf{DSN}$ of dynamical systems on networks of manifolds]
It is easy to see that if $\varphi:(G,{\mathcal P}) \to (G',{\mathcal P}')$ and $\psi:(G',{\mathcal P}')\to (G'',
{\mathcal P}'')$ are two f\/ibrations then
\begin{gather*}
(\psi\circ \varphi)^* = \varphi^*\circ \psi^*.
\end{gather*}
This can be interpreted as saying that the assignment
\begin{gather*}
(G,{\mathcal P})\mapsto \mathpzc{Ctrl}(G,{\mathcal P})
\end{gather*}
extends to a~contravariant functor $\mathpzc{Ctrl}$ from the category
$(\mathsf{Man}/\mathsf{Graph})_{\mbox{\sf{\tiny{f\/ib}}}} $ of networks of manifolds and f\/ibrations to the category
$\mathsf{Vect}$ of real vector spaces and linear maps.
That is, on arrows,
\begin{gather*}
\mathpzc{Ctrl}{\varphi}: = \varphi^*.
\end{gather*}
Grothendieck's construction (see for example~\cite{Awodey}) applied to this functor produces a~category $\mathsf{DSN}$
which we would like to call the category of (continuous time) dynamical systems on networks of manifolds.
More explicitly the objects of the category $\mathsf{DSN}$ are triples
\begin{gather*}
(G, {\mathcal P}\colon G_0 \to \mathsf{Man}, w\in \mathpzc{Ctrl}(G,{\mathcal P})),
\end{gather*}
where as before~$G$ is a~f\/inite directed graph, ${\mathcal P}$ is a~phase space function and $w=(w_a)_{a\in G_0}$ is
a~tuple of control systems associated with the input trees of the graph~$G$ and the function ${\mathcal P}$.

A morphism from $(G',{\mathcal P}', w')$ to $(G, {\mathcal P},w)$ is a~graph f\/ibration $\varphi\colon G\to G'$ with
${\mathcal P}' \circ \varphi = {\mathcal P}$ and $\varphi^*w' = w$.
Alternatively we may think of a~map from $(G',{\mathcal P}', w')$ to $(G, {\mathcal P},w)$ as a~f\/ibration of networks of
manifolds $\varphi\colon (G,{\mathcal P}) \to (G',{\mathcal P}')$ with $\varphi^*w' = w$.

Note that the Grothendieck construction also gives us a~forgetful functor
\begin{gather*}
\mathsf{DSN} \to \op{(\mathsf{Man}/\mathsf{Graph})_{\mbox{\sf{\tiny{f\/ib}}}}}
\end{gather*}
that simply forgets the open systems.
On objects it is given by sending the triple $(G,{\mathcal P}, w)$ to the pair $(G,{\mathcal P})$.
\end{Remark}

Of course just because we can def\/ine a~category and call it a~category of dynamical systems on networks does not mean
that this is a~right thing to do.
This said, Theorem~\ref{thm:3.8} tells us that to every dynamical system on a~network $(G, {\mathcal P}, w)$ we can
assign a~dynamical system $(\mathbb{P} G, \mathscr{I} w)$.
We will next argue that this assignment actually extends to a~functor
\begin{gather*}
\mathbb{P}: \ \mathsf{DSN} \to \mathsf{DS}
\end{gather*}
from dynamical systems on networks to the category $\mathsf{DS}$ of dynamical systems (q.v.\
Def\/inition~\ref{def:dyn-sys} and Remark~\ref{rmrk:dyn-sys}).
The f\/irst step is to def\/ine the functor on arrows.
We do it in Theorem~\ref{thm:main} below which may be considered the main result of the paper.

\begin{Theorem}
\label{thm:main}
Let $\varphi\colon (G,{\mathcal P})\to (G',{\mathcal P}')$ be a~fibration of networks of manifolds.
Then the pullback map
\begin{gather*}
\varphi^*\colon \ \mathpzc{Ctrl}(G',{\mathcal P}')\to \mathpzc{Ctrl}(G,{\mathcal P})
\end{gather*}
constructed in Theorem~{\rm \ref{thm:4.4}} is compatible with the interconnection maps
\begin{gather*}
\mathscr{I}'\colon \ \mathpzc{Ctrl}(G',{\mathcal P}') \to \Gamma (T\mathbb{P} G')
\qquad
\text{and}
\qquad
\mathscr{I}\colon \ \mathpzc{Ctrl}(G,{\mathcal P}) \to \Gamma (T\mathbb{P} G).
\end{gather*}
Namely for any collection $w'\in \mathpzc{Ctrl}(G',{\mathcal P}')$ of open systems on the network $(G', {\mathcal P}')$
the diagram
\begin{gather}
\label{eq:7}
\begin{split}
&
\xy
 (-15, 10)*+{T\mathbb{P} G'}="1";
 (15, 10)*+{T\mathbb{P} G}="2";
 (-15, -5)*+{\mathbb{P} G'}="3";
 (15, -5)*+{\mathbb{P} G}="4";
{\ar@{->}^{D\mathbb{P}\varphi } "1";"2"};
{\ar@{->}^{\mathscr{I}'(w') } "3";"1"};
{\ar@{->}_{\mathscr{I} (\varphi^*w') } "4";"2"};
{\ar@{->}_{\mathbb{P}\varphi} "3";"4"};
\endxy
\end{split}
\end{gather}
commutes.
Consequently
\begin{gather*}
\mathbb{P} \varphi\colon \ (\mathbb{P} (G',{\mathcal P}'), \mathscr{I}' (w'))\to (\mathbb{P} (G, {\mathcal P}),
\mathscr{I} (\varphi^* w'))
\end{gather*}
is a~map of dynamical systems.
\end{Theorem}

\begin{proof}
Recall that the manifold $\mathbb{P} G$ is the product $\bigsqcap_{a\in G_0} \mathbb{P} a$.
Hence the tangent bundle bundle $T\mathbb{P} G$ is the product $\bigsqcap_{a\in G_0} T\mathbb{P} a$.
The canonical projections
\begin{gather*}
T\mathbb{P} G\to T\mathbb{P} a
\end{gather*}
are the dif\/ferentials of the maps $\mathbb{P} \iota_a\colon \mathbb{P} G\to \mathbb{P} a$, where, as before,
$\iota_a\colon \{a\} \hookrightarrow G$ is the canonical inclusion of graphs.
Hence by the universal property of products, two maps into $T\mathbb{P} G$ are equal if and only if all their components
are equal.
Therefore, in order to prove that~\eqref{eq:7} commutes it is enough to show that
\begin{gather*}
D \mathbb{P} \iota_a \circ \mathscr{I} (\varphi^*w') \circ \mathbb{P} \varphi = D \mathbb{P} \iota_a \circ D \mathbb{P}
\varphi\circ \mathscr{I}'(w')
\end{gather*}
for all nodes $a\in G_0$.
By def\/inition of the restriction $\varphi|_{\{a\}}$ of $\varphi\colon G\to G'$ to $\{a\}\hookrightarrow G$, the diagram
\begin{gather}
\label{eq:**}
\begin{split}
& \xy
 (-10, 15)*+{\{a\}}="1";
 (10, 15)*+{\{\varphi(a)\}}="2";
 (-10, 0)*+{ G}="3";
 (10, 0)*+{ G'}="4";
{\ar@{->}^{\varphi|_{\{a\}} } "1";"2"};
{\ar@{->}_{\iota_a } "1";"3"};
{\ar@{->}_{\iota_{\varphi(a)} } "2";"4"};
{\ar@{->}_{\varphi} "3";"4"};
\endxy
\end{split}
\end{gather}
commutes.
By the def\/inition of the pullback map $\varphi^*$ and the interconnection maps $\mathscr{I}$, $\mathscr{I}'$ the diagram
\begin{gather}
\label{eq:8}
\begin{split}
&
\xy
 (-15, 15)*+{T\mathbb{P} a}="1";
 (15, 15)*+{T\mathbb{P}\varphi(a)}="2";
 (-15, -0)*+{\mathbb{P} I(a)}="3";
 (15, 0)*+{\mathbb{P} I(\varphi(a)) }="4";
 (-15, -15)*+{\mathbb{P}  G}="5";
(15, -15)*+{\mathbb{P}  G'}="6";
{\ar@{->}^{D\mathbb{P}\varphi|_{\{a\}} } "2";"1"};
{\ar@{->}_{(\varphi^*w')_a } "3";"1"};
{\ar@{->}^{w'_{\varphi(a)} } "4";"2"};
{\ar@{->}_{\mathbb{P}\xi_a } "5";"3"};
{\ar@{->}^{\mathbb{P}\xi_{\varphi(a)} } "6";"4"};
{\ar@{->}_{\mathbb{P}\varphi_a} "4";"3"};
{\ar@{->}^{\mathbb{P}\varphi} "6";"5"};
{\ar@/^{3pc}/^{\mathscr{I} (\varphi^*w')_a} "5";"1"};
{\ar@/_{3pc}/_{\mathscr{I}' (w')_{\varphi(a)}} "6";"2"};
\endxy
\end{split}
\end{gather}
commutes as well.
We now compute:
\begin{gather*}
D\mathbb{P}\iota_a \circ \mathscr{I} (\varphi^*w')\circ \mathbb{P}\varphi  =   (\mathscr{I} (\varphi^*w'))_a \circ
\mathbb{P}\varphi
\\
\qquad
  =   D \mathbb{P} (\varphi|_{\{a\}})\circ \mathscr{I}' (w')_{\varphi(a)}
\qquad
\hspace{21mm}
\text{by~\eqref{eq:8}}
\\
\qquad
  =   D \mathbb{P} (\varphi|_{\{a\}})\circ D\mathbb{P} \iota_{\varphi(a)} \circ \mathscr{I}' (w')
\hspace{18mm}
\text{by def\/inition of}~\mathscr{I}' (w')_{\varphi(a)}
\\
\qquad
  =   D \mathbb{P} \left(\iota_{\varphi(a)} \circ \varphi|_{\{a\}} \right) \circ \mathscr{I}' (w')
\hspace{23mm}
\text{since $\mathbb{P}$ is a~contravariant functor}
\\
\qquad
  =   D \mathbb{P} \left(\varphi \circ \iota_a \right) \circ \mathscr{I}' (w')
\hspace{34mm}
\text{by~\eqref{eq:**}}
\\
\qquad
  =   D \mathbb{P} (\iota_a) \circ D\mathbb{P} \varphi \circ \mathscr{I}' (w').
\end{gather*}
And we are done.
\end{proof}

\begin{Corollary}
The map
\begin{gather*}
\mathsf{DSN} \to \mathsf{DS}
\\
\big((G',{\mathcal P}', w') \xrightarrow{\varphi} (G,{\mathcal P},w)\big)
\ \mapsto \ \big((\mathbb{P} G', \mathscr{I} (w'))\xrightarrow{\mathbb{P} \varphi} (\mathbb{P} G, \mathscr{I} (w))\big)
\end{gather*}
is a~functor.
\end{Corollary}

\begin{Remark}
Given a~dynamical system on a~network $(G,{\mathcal P},w)$ we can forget the dynamics.
This def\/ines a~functor
\begin{gather*}
\mathsf{DSN} \to (\mathsf{Man}/\mathsf{Graph})_{\mbox{\sf{\tiny{f\/ib}}}}
\end{gather*}
from the category of dynamical systems on networks to a~subcategory of the category of networks of manifolds whose maps
are f\/ibrations of networks (hence the subscript $_{\mbox{\sf{\tiny{f\/ib}}}}$).
Composing the functor above with the functor $\mathsf{Man}/\mathsf{Graph} \to \op{\mathsf{Graph}}$ forgets all the
information except for the graph.
This gives rise to a~functor
\begin{gather*}
\mathsf{DSN}\to \opfib{\mathsf{Graph}}.
\end{gather*}
Here the superscript $\op{}$ indicates that the functor reverses the direction of arrows and the subscript $\gf{}$
reminds us that the morphisms in the target category are the (opposite of the) graph f\/ibrations.

These two functors from $\mathsf{DSN}$ to $\mathsf{DS}$ and to $\opfib{\mathsf{Graph}}$, respectively, allow us to
interpret continuous time dynamical systems on networks both as dynamical systems and as graphs.
\end{Remark}
We end the paper with examples.
\begin{Example}
\label{example:14}
Consider the graph f\/ibration
\begin{gather*}
\begin{tikzpicture}
[->,>=stealth',shorten >=1pt,auto,node distance=2cm, thick,main node/.style={circle,draw,font=\sffamily\bfseries}]
\node[main node] at (-2,0) (a1) {a$_1$}; \node[main node] at (-2,-2) (a2) {a$_2$}; \node[main node] at (0,-1) (bb) {b};
\node[main node] at (3,-1) (1) {a}; \node[main node] at (5,-1) (2) {b}; \node[main node] at (7,-1) (3) {c}; \node at
(1.5,-1) (arrow) {$\longrightarrow$};
\path[every node/.style={font=\sffamily\small}] (a1) edge node {$\gamma$} (bb) (a2) edge node [below] {$\delta$} (bb)
(1) edge [bend left] node {$\gamma'$} (2) (1) edge [bend right] node [below] {$\delta'$} (2) (2) edge node {}(3)
;
\end{tikzpicture}
\end{gather*}
as in Examples~\ref{example:5} and~\ref{example:6}.
Let ${\mathcal P}'\colon G_0'\to \mathsf{Man}$ be a~phase space function and let ${\mathcal P}= {\mathcal P}'\circ
\varphi$.
Then
\begin{gather*}
\mathbb{P} G' = {\mathcal P}'(a)\times {\mathcal P}' (b)\times {\mathcal P}' (c),
\\
\mathbb{P} G = {\mathcal P}'(a)\times {\mathcal P}' (a)\times {\mathcal P}' (b),
\\
\mathbb{P} \varphi (x,y,z) = (x,x,y),
\end{gather*}
and
\begin{gather*}
D \mathbb{P} \varphi (p,q, r) = (p,p,q).
\end{gather*}
For any $w' = (w_a', w_b', w_c') \in \mathpzc{Ctrl} (G', {\mathcal P}')$,
\begin{gather*}
\left(\mathscr{I}' (w')\right) (x,y,z) = (w'_a (x), w_b' (x,x, y), w_c' (y,z)),
\\
\varphi^* w' = (w_a', w_a', w_b'),
\\
\left(\mathscr{I} (\varphi^* w')\right) (x_1,x_2,y) = (w'_a (x_1), w'_a (x_2), w_b' (x_1,x_2, y))
\end{gather*}
and
\begin{gather*}
\left(\mathscr{I} (\varphi^* w')\circ \mathbb{P} \varphi\right) (x,y,z) = (w'_a (x), w'_a (x), w_b' (x,x, y))
\end{gather*}
while
\begin{gather*}
\left(D\mathbb{P} \varphi\circ \mathscr{I}' (w')\right) (x,y,z) = D\mathbb{P} \varphi (w'_a (x), w_b' (x,x, y), w'_c
(y,z)) = (w'_a (x), w'_a (x), w_b' (x,x, y)).
\end{gather*}
Hence
\begin{gather*}
\left(\mathscr{I} (\varphi^* w')\circ \mathbb{P} \varphi\right) = \left(D\mathbb{P} \varphi\circ \mathscr{I}'
(w')\right)
\end{gather*}
as expected.
\end{Example}

\begin{Example}
In Example~\ref{example:14} above the map $\varphi\colon G\to G'$ is neither injective nor surjective.
It can, of course, be factored as a~surjection $\psi\colon G\to G''$:
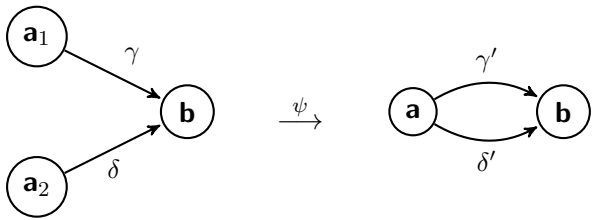
\begin{gather*}
\begin{tikzpicture}
[->,>=stealth',shorten >=1pt,auto,node distance=2cm, thick,main node/.style={circle,draw,font=\sffamily\bfseries}]
\node[main node] at (-2,0) (a1) {a$_1$}; \node[main node] at (-2,-2) (a2) {a$_2$}; \node[main node] at (0,-1) (bb) {b};
\node[main node] at (3,-1) (1) {a}; \node[main node] at (5,-1) (2) {b}; \node at (1.5,-1) (arrow)
{$\stackrel{\psi}{\longrightarrow}$};
\path[every node/.style={font=\sffamily\small}] (a1) edge node {$\gamma$} (bb) (a2) edge node [below] {$\delta$} (bb)
(1) edge [bend left] node {$\gamma'$} (2) (1) edge [bend right] node [below] {$\delta'$} (2)
;
\end{tikzpicture}
\end{gather*}
followed by an injection $\iota\colon G''\to G$:
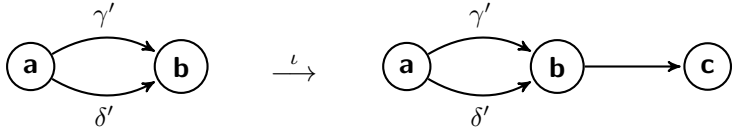
\begin{gather*}
\begin{tikzpicture}
[->,>=stealth',shorten >=1pt,auto,node distance=2cm, thick,main node/.style={circle,draw,font=\sffamily\bfseries}]
\node[main node] at (-2,-1) (a) {a}; \node[main node] at (0,-1) (bb) {b}; \node[main node] at (3,-1) (1) {a}; \node[main
node] at (5,-1) (2) {b}; \node[main node] at (7,-1) (3) {c}; \node at (1.5,-1) (arrow)
{$\stackrel{\iota}{\longrightarrow}$};
\path[every node/.style={font=\sffamily\small}] (a) edge [bend left] node {$\gamma'$} (bb) (a) edge [bend right] node
[below] {$\delta'$} (bb) (1) edge [bend left] node {$\gamma'$} (2) (1) edge [bend right] node [below] {$\delta'$} (2)
(2) edge node {}(3)
;
\end{tikzpicture}
\end{gather*}

The map $\mathbb{P} \psi\colon \mathbb{P} G'' \to \mathbb{P} G$ is easily seen to be given~by
\begin{gather*}
\mathbb{P}\psi (x,y) = (x,x,y).
\end{gather*}
It is an embedding, as it should be (q.v.\ Lemma~\ref{lemma:2.surj}).
The map $\mathbb{P} \imath\colon \mathbb{P} G'\to \mathbb{P} G''$ is given~by
\begin{gather*}
\mathbb{P} \imath (x,y,z) = (x,y).
\end{gather*}
It is a~submersion (q.v.\ Lemma~\ref{lemma:2.inj}).
Since $\mathbb{P}$ is a~contravariant functor,
\begin{gather*}
\mathbb{P} \varphi = \mathbb{P} (\imath \circ \psi) = \mathbb{P}\psi \circ \mathbb{P} \imath.
\end{gather*}
Theorem~\ref{thm:main} tells us that for any $w' = (w_a', w_b', w_c') \in \mathpzc{Ctrl} (G', {\mathcal P}')$, the map
$\mathbb{P} \imath$ projects the integral curves of the vector f\/ield $\mathscr{I} (w')$ to the integral curves of the
vector f\/ield $\mathscr{I} (\imath^* w')$ on~$\mathbb{P} G''$.
Furthermore, $\mathbb{P} \psi$ embeds the dynamical system $(\mathbb{P} G'', \mathscr{I} (\imath^* w'))$ into the
dynamical system $(\mathbb{P} G, \mathscr{I} (\varphi^* w'))$.
An interested reader can check these two assertions directly.
\end{Example}

\begin{Example}
Consider the injective graph f\/ibration $\iota\colon G\to G'$:
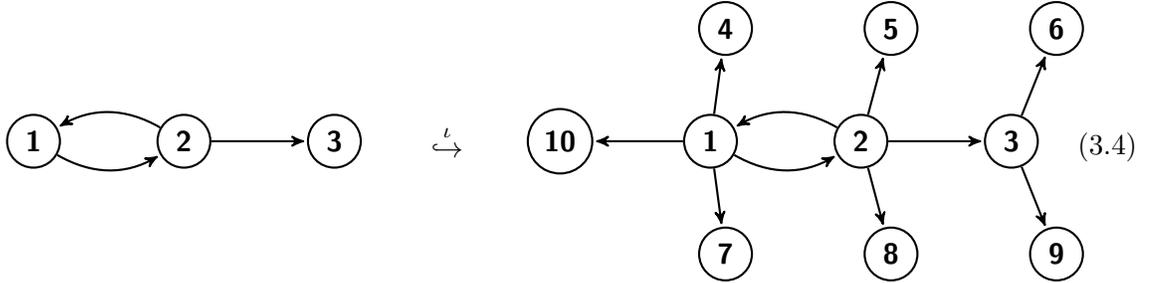
\begin{gather}
\begin{split}
& \begin{tikzpicture}[->,>=stealth',shorten >=1pt,auto,node distance=2cm,
  thick,main node/.style={circle,draw,font=\sffamily\bfseries}]
  \node[main node] at (2,-1) (1) {1};
  \node[main node] at (4,-1) (2) {2};
  \node[main node] at (6,-1) (3) {3};
  \node at (7.5,-1) (i) {$\stackrel{\iota}{\hookrightarrow}$};
  \node[main node] at (9,-1) (10') {10};
  \node[main node] at (11,-1) (1') {1};
  \node[main node] at (13,-1) (2') {2};
  \node[main node] at (15,-1) (3') {3};
  \node[main node] at (11.2,0.5) (4') {4};
  \node[main node] at (13.4,0.5) (5') {5};
  \node[main node] at (15.6,0.5) (6') {6};
  \node[main node] at (11.2,-2.5) (7') {7};
  \node[main node] at (13.4,-2.5) (8') {8};
  \node[main node] at (15.6,-2.5) (9') {9};
  \path[every node/.style={font=\sffamily\small}]
  (2) edge [bend right] node {} (1)
  (1) edge [bend right] node [below] {} (2)
  (2) edge  node {}(3)
  (2') edge [bend right] node {} (1')
  (1') edge [bend right] node [below] {} (2')
  (2') edge  node {}(3')
  (1') edge  node {}(10')
  (1') edge  node {}(4')
  (1') edge  node {}(7')
  (2') edge  node {}(5')
  (2') edge  node {}(8')
  (3') edge  node {}(6')
  (3') edge  node {}(9');
\end{tikzpicture}
\end{split}\!\!\!\!\!\label{eq:ten}
\end{gather}

Choose phase space functions ${\mathcal P}$, ${\mathcal P}'$ so that $i\colon (G,{\mathcal P})\to (G',{\mathcal P}')$ is
a~map of networks.
By Theorem~\ref{thm:main}, for any collection $w'\in \mathpzc{Ctrl}(G',{\mathcal P}')$ of open systems on the network
$(G', {\mathcal P}')$ the dynamics in the subsystem $(\mathbb{P} G, \mathscr{I}(i^*w'))$ drives the entire system
$(\mathbb{P} G', \mathscr{I} (w'))$.
This is intuitively clear from the graph~\eqref{eq:ten} since there are no ``feedbacks'' from vertices $4, \ldots,10$
back into~$1$,~$2$,~$3$.
\end{Example}

\subsection*{Acknowledgments}

L.D.\ was supported by the National Science Foundation under grants CMG-0934491 and UBM-1129198 and by the National
Aeronautics and Space Administration under grant NASA-NNA13 AA91A.
The authors also thank the anonymous referees whose comments signif\/icantly improved the manuscript.

\pdfbookmark[1]{References}{ref}
\LastPageEnding

\end{document}